\documentclass[12pt, reqno]{amsart}
\usepackage{amsmath}
\usepackage{amsthm}
\usepackage{amssymb}
\usepackage{enumerate}
\usepackage{graphicx}
\usepackage{graphicx}
\usepackage{amssymb,amsfonts}
\usepackage{amsmath}
\usepackage[colorlinks,linkcolor=blue,anchorcolor=blue,citecolor=blue]{hyperref}
\usepackage[numbers,sort&compress]{natbib}
\usepackage{marvosym}
\usepackage{epstopdf}
\usepackage{caption}
\usepackage{multicol}
\usepackage{multirow}
\usepackage{makebox}
\usepackage{subfigure}
\numberwithin{equation}{section}
\newtheorem{theo}{Theorem} 

\newtheorem{lem}{Lemma}

\newtheorem{defn}{Definition}
\begin{document}

\title{$\textrm{SL}(n)$ Contravariant Matrix-Valued Valuations on Polytopes}
\author[Zeng and Zhou]{Chunna Zeng, Yuqi Zhou}

\address{1.School of Mathematics  Sciences,
 Chongqing Normal University, Chongqing 401331, People's Republic of China; Institut f\"{u}r Diskrete Mathematik und Geometrie, Technische Universit\"{a}t Wien, Wiedner Hauptstra\ss e 8--10/1046, 1040 Wien, Austria }
\email{zengchn@cqnu.edu.cn}

\address{2.School of Mathematical Sciences, Chongqing Normal University, Chongqing 401332, People's Republic of China}
\email{zhouyuqi202212@163.com}

\thanks{Supported in part by the Major Special Project of NSFC (Grant No. 12141101), the Young Top-Talent program of Chongqing (Grant No. CQYC2021059145), the Research Project of Chongqing Education Commission CXQT21014 and the Characteristic innovation projects of universities in Guangdong province (Grant No. 2020KTSCX358)}
\thanks{{*}Corresponding author: Chunna Zeng}
\begin{abstract}
All $\textrm{SL}(n)$ contravariant matrix-valued valuations on polytopes in $\mathbb{R}^n$ are completely classified without any continuity assumptions. Moreover, the symmetry assumption of matrices is removed. The general Lutwak-Yang-Zhang matrix turns out to be the only such valuation if $n\geq 4$, while a new function shows up in dimension three. In dimension two, the classification corresponds to the known case of $\textrm{SL}(2)$ equivariant matrix-valued valuations.

\end{abstract}

\maketitle

\section{Introduction }

\par A classical concept from mechanics is the Legendre ellipsoid $\Gamma_2 K$ associated with a convex body $K$. The Legendre ellipsoid  is the unique ellipsoid
centered at the center of mass of $K$ such that the ellipsoid's moment of inertia about any axis passing through the center of mass is the same as that of $K$.
Write a vector $x \in \mathbb{R}^n$ as $x=(x_1,x_2,\ldots,x_n)$. The Legendre ellipsoid $\Gamma_2 K$ can be defined by the moment matrix $M(K)$ of $K$, that is, the $(n\times n)$-matrix with coefficients
\begin{align*}
M_{ij}(K)=\int_K x_i x_j dx.
\end{align*}
For any convex body $K$ with nonempty interior, its moment matrix $M(K)$ is positive definite.
The moment matrix is closely related to the Bourgain slicing problem,
which is one of  open problems in the asymptotic theory. The Bourgain slicing problem mainly
seeks a unique upper bound estimate for the isotropic constant. Klartag
showed that an upper bound of the isoropic constant has the dimension dependency $n^{\frac{1}{4}}$ (see \cite{Artstein,Klartag-Milman}).
 A recent advancement over Klartag's finding was contributed by Chen \cite{Chen}.

\par In \cite{L-Y-Z2}, Lutwak, Yang and Zhang defined a new ellipsoid $\Gamma_{-2}K$ for $K \in \mathcal{K}_{(o)}^n$, the space of convex bodies containing the origin in their interiors. Actually, the Lutwak-Yang-Zhang $(\textrm{LYZ})$ ellipsoid  is recognized as the dual counterpart of the Legendre ellipsoid in the Brunn-Minkowski theory. They \cite{L-Y-Z1} also established that $\Gamma_{-2}K \subset \Gamma_{2}K$ and noted that this was  a geometric analogue of the Cramer-Rao inequality. For a convex polytope $P \in \mathcal{P}_{(o)}^n,$  the space of convex polytopes containing the origin in their interiors in $\mathbb{R}^n.$ Then $\Gamma_{-2}P$ can be defined by the $\textrm{LYZ}$ matrix with coefficients (see \cite{L-Y-Z2})
\begin{align*}
L_{ij}(P)=\sum_{u\in \mathcal{N}(P)} \frac{a_P(u)}{h_P(u)} u_i u_j,
\end{align*}
where $\mathcal{N}(P)$ denotes the set of all outer unit normals of facets of $P$. For a unit normal $u\in \mathbb {S}^{n-1}$, denote by $a_P(u)$  the $(n-1)$-dimensional volume of the facet with normal $u$, and $h_P(u)=\max \{x\cdot u:x \in P\}$ the support function of $P$.
For additional information on the $\textrm{LYZ}$ ellipsoid and its connection to the Fisher information from information theory, see \cite{G-L-Y-Z, Ludwig4,L-Y-Z1,L-Y-Z2}.
Beyond $\Gamma_{2}$ and $\Gamma_{-2}$, there exist other well-known ellipsoids
like the John ellipsoid, the $L_{p}$ John ellipsoid, the Petty ellipsoid, and the $M$-ellipsoid. However, only $\Gamma_{2}$ and $\Gamma_{-2}$ are linear, and their  corresponding  moment matrix and LYZ matrix function are `valuations'.

\par A function $\mu: \mathcal{S}\rightarrow \langle\mathcal{A},+\rangle$ defined on a collection $\mathcal{S}$ of sets and taking values in an abelian semigroup $\langle\mathcal{A},+\rangle$ is called a valuation if
\begin{equation*}
\mu(P)+\mu(Q)=\mu(P\cup Q)+\mu(P\cap Q),
\end{equation*}
whenever $P,\,Q,\,P\cap Q,$ and $P\cup Q\in \mathcal{S}$.
At the beginning of the twentieth century, valuations were first constructed by
Dehn in his solution of Hilbert's Third Problem. Nearly 50 years later, Hadwiger
initiated a systematic study of valuations by his celebrated characterization theorem. He  revealed that all continuous, rigid motion invariant valuations on the
space of convex bodies (i.e., compact convex sets) in $\mathbb {R}^n$ are linear combinations of
intrinsic volumes. Valuations have various important applications in integral geometry (see \cite[Chap. 7]{Gruber}, \cite{Klain-Rota}, \cite[Chap. 6]{Ludwig3}). They also turned out to be extremely
fruitful and useful especially in the affine geometry of convex bodies. Examples of valuations are affine surface areas, intrinsic volumes, the intersection bodies, the projection bodies, and other Minkowski valuations. See also \cite{Alesker1,Alesker2,Alesker3,Alesker4,Klain1,Klain2,Klain3,Klain-Rota,Ludwig5,McMullen1,McMullen2,Zeng-Ma} for a historical account and some of the more recent contributions.

\par Denote by $\mathbb{M}_e^n$ the set of real symmetric $(n \times n)$-matrices, and  $\mathbb{M}^n$ the set of real $(n \times n)$-matrices over $\mathbb{R}^n.$
 A function $\mu:\mathcal{P}_{(o)}^n\rightarrow \mathbb{M}_e^n$ is called $\textrm{GL}(n)$ contravariant if there exists a constant $q \in \mathbb{R}$ such that
\begin{equation*}
\mu(\phi p)=|\det \phi^{-t}|^q \,\phi^{-t} \,\mu(p)\,\phi^{-1}
\end{equation*}
for all $P \in \mathcal{P}_{(o)}^n$ and $\phi \in \textrm{GL}(n)$. Let $\mathcal{P}^n$ be the set of  convex ploytopes, and $\mathcal{P}_o^n$ the subspace of  all convex ploytopes  containing the origin in $\mathbb{R}^n$.
 Denote by  $\mathcal{Q}^n$ be either $\mathcal{P}^n$ or $\mathcal{P}_o^n.$

\par A function $\mu:\mathcal{Q}^n\rightarrow \mathbb{M}^n$ is called $\textrm{SL}(n)$ contravariant if
\begin{align*}
\mu(\phi p)=\phi^{-t} \,\mu(p)\,\phi^{-1}
\end{align*}
for all $P \in \mathcal{Q}^n$ and $\phi \in \textrm{SL}(n)$, and is called $\textrm{SL}(n)$ equivariant if
\begin{align*}
\mu(\phi p)=\phi\,\mu(p)\,\phi^t
\end{align*}
for all $P \in \mathcal{Q}^n$ and $\phi \in \textrm{SL}(n)$.

A function defined on $\mathcal{P}_{(o)}^n$ is called (Borel) measurable if the preimage of every open set is a Borel set.
In 2003, Ludwig \cite{Ludwig4} established the first characterization of the moment matrix and  $\textrm{LYZ}$ matrix.
\begin{lem}\cite{Ludwig4}\label{lem0}
Let $n\geq 3$. A function $\mu: \mathcal{P}_{(o)}^n \rightarrow \mathbb{M}_e^n$ is a measurable $\textrm{GL}(n)$ contravariant valuation if and only if there exists a constant $c\in \mathbb{R}$ such that
\begin{equation*}
\mu(P)=c\,M(P^*) \qquad \textrm{or} \qquad \mu(P)=c\,L(P)
\end{equation*}
for every $P \in \mathcal{P}_{(o)}^n$, where $P^*$ is the polar body of $P$.
\end{lem}

\par Recently, Ma and Wang \cite{Ma-Wang} extended the $\textrm{LYZ}$ matrix from convex ploytopes containing the origin in their interiors to arbitrary convex ploytopes. For a solution of Cauchy's functional equation $\zeta:\mathbb{R}\rightarrow \mathbb{R}$, the general $\textrm{LYZ}$ matrix $L_\zeta(P)$ of $P \in \mathcal{P}^n$ is the $(n\times n)$-matrix with coefficients
\begin{equation*}
L_{\zeta,ij}(P)=\sum_{u\in \mathcal{N}(P) \setminus \{h_P=0\}} \frac{\zeta(a_P(u)h_P(u))}{h_P^2 (u)} u_i u_j.
\end{equation*}
The general LYZ matrix is also shown to be essentially the unique class of $\textrm{SL}(n)$ contravariant symmetric matrix-valued valuations. Furthermore, they extended Ludwig's result to $\mathcal{P}_o^n$ without any homogeneity assumptions or any continuity assumptions.
\begin{lem}\cite{Ma-Wang}\label{lemLYZ}
Let $n\geq 3$. A function $\mu: \mathcal{P}_o^n \rightarrow \mathbb{M}_e^n$ is an $\textrm{SL}(n)$ contravariant valuation if and only if there exists a solution of Cauchy's functional equation $\zeta:\mathbb{R} \rightarrow \mathbb{R}$ such that
\begin{equation*}
\mu(P)=L_\zeta (P)
\end{equation*}
for every $P \in \mathcal{P}_o^n$.
\end{lem}

\par The aim of this paper is to obtain a complete classification of $\textrm{SL}(n)$ contravariant matrix-valued valuations on polytopes. Moreover, the symmetry assumption of matrices is removed.

\begin{theo}\label{theo1}
Let $n\geq 4$. A function $\mu: \mathcal{P}_o^n \rightarrow \mathbb{M}^n$ is an $\textrm{SL}(n)$ contravariant valuation if and only if there exists a solution of Cauchy's functional equation $\zeta:\mathbb{R} \rightarrow \mathbb{R}$ such that
\begin{equation*}
\mu(P)=L_\zeta (P)
\end{equation*}
for every $P \in \mathcal{P}_o^n$.
\end{theo}

\par Let $(\textbf{i},\,\textbf{j},\,\textbf{k})$ be a positively oriented orthonormal basis. Each vector $\textbf{a}\in\mathbb{R}^3$ can be define as the sum of three orthogonal components parallel to the standard basis vectors, that is,
\begin{align*}
\textbf{a}=a_1\textbf{i}+a_2\textbf{j}+a_3\textbf{k}.
\end{align*}
Then the antisymmetric matrix $\textrm{Anti}(\textbf{a})$, for a vector $\textbf{a}$, is defined by
\begin{align*}
\textrm{Anti}(\textbf{a})= \begin{pmatrix}
0 & -a_3 & a_2 \\
a_3 & 0 & -a_1 \\
-a_2 & a_1 & 0
\end{pmatrix}.
\end{align*}
\par For $\textbf{a}=a_1\textbf{i}+a_2\textbf{j}+a_3\textbf{k}$ and $\textbf{b}=b_1\textbf{i}+b_2\textbf{j}+b_3\textbf{k}\in \mathbb{R}^3$, the cross product $\textbf{a} \times \textbf{b}$ is defined by
\begin{align*}
\textbf{a} \times \textbf{b}=
\begin{vmatrix}
\textbf{i} & \textbf{j} & \textbf{k}\\
a_1 & a_2 & a_3\\
b_1 & b_2 & b_3
\end{vmatrix}.
\end{align*}
It also can  be expressed as the product of an antisymmetric matrix and a vector, i.e.,
\begin{align*}
\textbf{a}\times \textbf{b}=\textrm{Anti}(\textbf{a}) \textbf{b}= \begin{pmatrix}
0 & -a_3 & a_2 \\
a_3 & 0 & -a_1 \\
-a_2 & a_1 & 0
\end{pmatrix}\begin{pmatrix}
b_1 \\
b_2 \\
b_3
\end{pmatrix}.
\end{align*}
More generally, it obeys the following identity under matrix transformation $\phi \in \textrm{SL}(3)$
\begin{align*}
(\phi \textbf{a}) \times (\phi \textbf{b})=\phi^{-t}(\textbf{a}\times \textbf{b}).
\end{align*}
Thus, we obtain
\begin{align*}
\textrm{Anti}(\phi \textbf{a})=\phi^{-t} \textrm{Anti}(\textbf{a})\phi^{-1},
\end{align*}
which implies that $\textrm{Anti}$ is $\textrm{SL}(3)$ contravariant.
\par Let $\mathcal{T}^3$ be the set of simplices in $\mathbb{R}^3$ with one vertex at the origin. Define  function $I:\mathcal{T}^3 \rightarrow \mathbb{M}^3$ as
\begin{align*}
I(T)=\det(u,v,w) \textrm{Anti}\left((u,v,w)\textbf{e}\right)
\end{align*}
if $T \in \mathcal{T}^3$ with dim $T=3$ and $T=[o,u,v,w]$ with $\det(u,v,w) >0$, where $\textbf{e}=\textbf{i}+\textbf{j}+\textbf{k}$; or
\begin{align*}
I(T)=\textbf{0}
\end{align*}
if $T \in \mathcal{T}^3$ with dim $T\leq 2$. Here $\textbf{0}$ denotes the matrix where every element is $0$. In Section 2, it will be derived that $I$ is an $\textrm{SL}(3)$ contravariant valuation.

\begin{theo}\label{theo2}
A function $\mu: \mathcal{P}_o^3 \rightarrow \mathbb{M}^3$ is an $\textrm{SL}(3)$ contravariant valuation if and only if there exist a constant $c \in \mathbb{R}$ and a solution of Cauchy's functional equation $\xi: [0, \infty)\rightarrow \mathbb{R}$ such that
\begin{align}\label{3.2.1}
\mu(P)=c\,I(P)+2\,L_\xi(P)
\end{align}
for every $P\in \mathcal{P}_o^3$. The notation $I(P)$ stands for $I(P)=\sum\limits_{i=1}^{m}I(T_i)$, where $\{T_1,\ldots, T_m\}$ is a triangulation of $P$ that is dissected at the origin.
\end{theo}

\par Using a relation with $\textrm{SL}(2)$ equivariant matrix-valued valuation, we obtain the classification in the two-dimensional case.

\begin{theo}\label{theo3}
A function $\mu: \mathcal{P}_o^2 \rightarrow \mathbb{M}^2$ is an $\textrm{SL}(2)$ contravariant valuation if and only if there exist constants $c_1, c_2, c_3, c_4 \in \mathbb{R}$ and solutions of Cauchy's functional equation $\alpha, \beta: [0, \infty)\rightarrow \mathbb{R}$ such that
\begin{align}
\mu(P)=c_1\,\tilde{M}(P)+c_2\,\tilde{E}(P)+c_3\,\tilde{F}(P)+\tilde{H}_\alpha(P)+\tilde{G}_\beta(P)+c_4\,\rho_{\pi/ 2}
\end{align}
for every $P\in \mathcal{P}_o^2 $, where $\rho_{\pi/ 2}$ is the counterclockwise rotation in $\mathbb{R}^2$ by the angle $\pi/ 2$.
\end{theo}

\par Here, $\tilde{M}(P)=\psi_{\pi/ 2} \,M(P)\, \psi_{\pi/ 2}^{-1}$, where $\psi_{\pi/ 2} \in \mathbb{M}^2$ denotes the rotation by the angle $\pi/ 2$. The functions $E(P), F(P), H_\alpha(P)$ and  $G_\beta(P)$ are $\textrm{SL}(2)$ equivariant valuations which were proved in \cite{Wang}.  By applying the same transform to $E(P),\,F(P),\,H_\alpha(P)$, and $G_\beta(P)$, we obtain $\tilde{E}(P),\,\tilde{F}(P),\,\tilde{H}_\alpha(P)$, and $\tilde{G}_\beta(P)$, respectively.

\par Next, we  extend these results to $\mathcal{P}^n$. This step is as in the classification of convex body valuations by Schuster and Wannerer \cite{Schuster} and Wannerer \cite{Wannerer}.

\begin{theo}\label{theo4}
Let $n\geq 4$. A function $\mu: \mathcal{P}^n \rightarrow \mathbb{M}^n$ is an $\textrm{SL}(n)$ contravariant valuation if and only if there exist solutions of Cauchy's functional equation $\zeta_1, \zeta_2 :\mathbb{R} \rightarrow \mathbb{R}$ such that
\begin{equation*}
\mu(P)=L_{\zeta_1} (P)+L_{\zeta_2} ([o, P])
\end{equation*}
for every $P \in \mathcal{P}^n$, where $[o, P]$ denotes the convex hull of the origin and $P$.
\end{theo}

\par Again, the cases of dimension three and dimension two are different. We show the following results.

\begin{theo}\label{theo5}
If $\mu: \mathcal{P}^3 \rightarrow \mathbb{M}^3$ is an $\textrm{SL}(3)$ contravariant valuation if and only if there exist constants $c_1, c_2 \in \mathbb{R}$ and solutions of Cauchy's functional equation $\xi_1, \xi_2: [0, \infty)\rightarrow \mathbb{R}$ such that
\begin{align}\label{1.3}
\mu(P)=c_1\, I([o,P])+c_2\, \sum_{i=1}^{l}I([o,F_i])+2\, L_{\xi_1}([o,P])+2\, \sum_{i=1}^{l} L_{\xi_2}([o,F_i])
\end{align}
for every $P \in \mathcal{P}^3$ with facets $F_1,\ldots, F_l$ visible from the origin.
\end{theo}

\begin{theo}\label{theo6}
A function $\mu: \mathcal{P}^2 \rightarrow \mathbb{M}^2$ is an $\textrm{SL}(2)$ contravariant valuation if and only if there exist constants $c_1, c_2, c_3, c_4, c_5 \in \mathbb{R}$ and solutions of Cauchy's functional equation $\alpha, \alpha', \beta, \beta': [0, \infty)\rightarrow \mathbb{R}$ such that
\begin{align}
\mu(P)=&c_1\,\tilde{M}(P)+c_2\,\tilde{M}([o,P])+c_3\,\tilde{E}([o,P])+c_4\,\tilde{F}([o,P])+
\tilde{H}_\alpha([o,P])\notag\\
&+\sum_{i=1}^{r-1} \tilde{H}_{\alpha'}([o,v_{i+1},v_i])+ \tilde{G}_\beta([o,P])+\sum_{i=1}^{r-1}\tilde{G}_{\beta'}([o,v_{i+1},v_i])+c_5\,\rho_{\pi/2}
\end{align}
for every $P\in \mathcal{P}^2 $ with vertices $v_1, \ldots , v_r$ visible from the origin and labeled counterclockwise.
\end{theo}

\section{Notation and Preliminaries}
\par We work in $n$-dimensional Euclidean space $\mathbb{R}^n$ with origin $o$, and denote its standard basis by $e_1,\ldots,e_n$. Write the coordinates of a vector $x \in \mathbb{R}^n$ with respect to the standard basis by $x_1,x_2,\ldots,x_n$. The standard inner product of $x, y \in \mathbb{R}^n$ is denoted by $x\, \cdot\, y$. Denote the determinant of a matrix $A$ by $\det A$, and the $n\times n$ identity matrix by $I_n=(e_1,\ldots,e_n)$. The affine hull, the relative interior, the interior, the dimension, and the boundary of a given set in $\mathbb{R}^n$ are denoted by aff, relint, int, dim, and bd, respectively.
\par Denote the convex hull of $v_1, v_2,\ldots,v_k \in \mathbb{R}^n$ by $[v_1, v_2,\ldots,v_k]$. A convex polytope is the convex hull of finitely many points in $\mathbb{R}^n$. Two basic classes of polytopes are the $k$-dimensional standard simplex $T^k=[o,e_1,\ldots,e_k]$ and $\tilde{T}^{k}=[e_1,\ldots,e_k]$, which is a $(k-1)$-dimensional simplex. For $ i=1,2,\ldots,n$, let $\mathcal{T}^i$ be the set of $k$-dimensional simplices with one vertex at the origin, and let $\tilde{\mathcal{T}}^{i-1}$ be the set of $(i-1)$-dimensional simplices $T \in \mathbb{R}^n$ with $0\not\in$  $T$. Generally, there exists a triangulation of a $k$-dimensional polytope $P$ into simplices as a set of $k$-dimensional simplices $\{T_1,\ldots,T_r\}$ which have pairwise disjoint interiors, with $P = \bigcup_i T_i$ and with the property that, for an arbitrary $1 \leq i_1 < \cdots <i_j \leq r$, the intersections $T_{i_1} \cap \cdots \cap T_{i_j}$ are again simplices.

\par For more basic results on valuations, see \cite{Klain-Rota} and \cite{Parapatits}. We now recall some basic results on valuations.

\begin{lem}\cite{Parapatits}\label{lem 2}
Let $\mu: \mathcal{P}_o^n \rightarrow \mathcal{A}$ be a valuation. Then $\mu$ is determined by its values on $n$-simplices with one vertex at the origin and its value on $\{o\}$.
\end{lem}

\begin{lem}\cite{Parapatits}\label{lem 3}
Let $\mu: \mathcal{P}^n \rightarrow \mathcal{A}$ be a valuation. Then $\mu$ is determined by its values on $n$-simplices.
\end{lem}

\par  Denote by $\textrm{SL}(n)$ the group of special linear transformations, i.e., those with linear transformations $\phi$ with determinant $\det \phi=1$, and by $\textrm{GL}(n)$ the group of general linear transformations, i.e., those with linear transformations $\phi$ with determinant $\det \phi \neq0$.  In addition, a general reference on convex geometry is the book by Schneider \cite{Schneider} or Gardner \cite{Gardner1}. Let $\mu:\mathcal{Q}^n\rightarrow \mathbb{M}^n$ be an $\textrm{SL}(n)$ contravariant valuation. Then $\mu$ can be decomposed as
\begin{equation*}
\mu(P)=\frac{1}{2}(\mu(P)+\mu(P)^t)+\frac{1}{2}(\mu(P)-\mu(P)^t)
\end{equation*}
for all $P \in \mathcal{Q}^n$. Clearly, $\mu(P)+\mu(P)^t$ is a symmetric matrix and $\mu(P)-\mu(P)^t$ is an antisymmetric matrix.

\par A valuation on $\mathcal{Q}^n$ is called simple if it vanishes on every $P \in \mathcal{Q}^n$ with dim $P \leq n-1$.

\begin{lem}\label{lem2.1}
The function $I$ is a simple $\textrm{SL}(3)$ contravariant valuation on $\mathcal{T}^3$.
\end{lem}

\emph{Proof.} \
By the definition of $I$ we obtain that it is indeed  simple on $\mathcal{T}^3$.
\par Next, we are going to prove that $I$ is a valuation on $\mathcal{T}^3$. It suffices to show that
\begin{align}\label{2-3}
I(T_1)+I(T_2)=I(T_1 \cup T_2)+I(T_1 \cap T_2)
\end{align}
for $T_1, T_2 \in \mathcal{T}^3$ with $ T_1 \cup T_2 \in \mathcal{T}^3$. It is clear that (\ref{2-3}) holds for $T_1 \subseteq T_2$ or $T_1 \supseteq T_2$. So we just need to consider  the case dim $T_1=$ dim $T_2=3$ with $T_1 \nsubseteq T_2$ and $T_1 \nsupseteq T_2$.
\par Assume that 
\begin{align*}
T_1=[o,u,v,k_1v+(1-k_1)w]
\end{align*}
and
\begin{align*}
T_2=[o,u,k_2v+(1-k_2)w,w],
\end{align*}
where $0 < k_1 \leq k_2 <1$ and $\det(u,v,w)>0$. Without loss
of generality, we have
\begin{align*}
T_1 \cup T_2=[o,u,v,w]
\end{align*}
and
\begin{align*}
T_1 \cap T_2=[o,u,k_2v+(1-k_2)w,k_1v+(1-k_1)w].
\end{align*}
Setting $\det(u,v,w)=r$. We have $\det(u,v,k_1v+(1-k_1)w)=(1-k_1) r$, $\det(u,k_2v+(1-k_2)w,w)=k_2r$, and $\det(u,k_2v+(1-k_2)w,k_1v+(1-k_1)w)=(k_2-k_1) r$. Therefore,
\begin{align*}
&I(T_1)=(1-k_1) r \, \textrm{Anti}\left((u,v,k_1v+(1-k_1)w)\textbf{e}\right),\\
&I(T_2)=k_2r  \, \textrm{Anti}\left((u,k_2v+(1-k_2)w,w)\textbf{e}\right),\\
&I(T_1\cup T_2)=r \,\textrm{Anti}\left((u,v,w)\textbf{e}\right),
\end{align*}
and
\begin{align*}
I(T_1\cap T_2)=(k_2-k_1) r \, \textrm{Anti}\left((u,k_2v+(1-k_2)w,k_1v+(1-k_1)w)\textbf{e}\right).
\end{align*}
Let $u=(u_1, u_2, u_3)^t$, $v=(v_1, v_2, v_3)^t$ and $w=(w_1, w_2, w_3)^t$, then
\begin{align}\label{2.4}
&I(T_1)+I(T_2)\\ \nonumber
=&(1-k_1)r \,\textrm{Anti}\left(\begin{array}{cl}
u_1+(1+k_1)v_1 +(1-k_1)w_1\\
u_2+(1+k_1)v_2 +(1-k_1)w_2\\
u_3+(1+k_1)v_3 +(1-k_1)w_3\\
\end{array}\right)\\ \nonumber
&\ \ \ \ +k_2r \,\textrm{Anti}\left(\begin{array}{cl}
u_1+k_2v_1 +(2-k_2)w_1\\
u_2+k_2v_2 +(2-k_2)w_2\\
u_3+k_2v_3 +(2-k_2)w_3\\
\end{array}\right)\\ \nonumber
=&\textrm{Anti}\left(\begin{array}{cl}
\left(1-k_1+k_2\right)ru_1 + \left(1-k_1^2+k_2^2\right)rv_1 + ((1-k_1)^2+(2k_2-k_2^2))rw_1 \\
\left(1-k_1+k_2\right)ru_2 + \left(1-k_1^2+k_2^2\right)rv_2 + ((1-k_1)^2+(2k_2-k_2^2))rw_2 \\
\left(1-k_1+k_2\right)ru_3 + (\left(1-k_1^2+k_2^2\right)rv_3 + ((1-k_1)^2+(2k_2-k_2^2))rw_3
\end{array}\right),\nonumber
\end{align}
and
\begin{align}\label{2.5}
&I(T_1\cup T_2)+I(T_1\cap T_2)\\ \nonumber
=&r \,\textrm{Anti}\left(\begin{array}{cl}
u_1+v_1 +w_1\\
u_2+v_2 +w_2\\
u_3+v_3 +w_3\\
\end{array}\right)+(k_2-k_1)r \,\textrm{Anti}\left(\begin{array}{cl}
u_1+(k_1+k_2)v_1 +(2-k_1-k_2)w_1\\
u_2+(k_1+k_2)v_2 +(2-k_1-k_2)w_2\\
u_3+(k_1+k_2)v_3 +(2-k_1-k_2)w_3\\
\end{array}\right)\\ \nonumber
=&\textrm{Anti}\left(\begin{array}{cl}
\left(k_2-k_1+1\right)ru_1 + \left(k_2^2-k_1^2+1\right)rv_1 + ((k_2-k_1)(2-k_1-k_2)+1)rw_1 \\
\left(k_2-k_1+1\right)ru_2 + \left(k_2^2-k_1^2+1\right)rv_2 + ((k_2-k_1)(2-k_1-k_2)+1)rw_2 \\
\left(k_2-k_1+1\right)ru_3 + (\left(k_2^2-k_1^2+1\right)rv_3 + ((k_2-k_1)(2-k_1-k_2)+1)rw_3
\end{array}\right).
\end{align}
Comparing coefficients in  (\ref{2.4}) and (\ref{2.5})  shows that  (\ref{2-3}) holds.
\par Finally, we show that $I$ is $\textrm{SL}(3)$ contravariant. Let $T \in \mathcal{T}^3$ with dim $T=3$ and $T=[o,u,v,w]$ with $\det(u,v,w) >0$. Then for $\phi \in \textrm{SL}(3)$, we have
\begin{align*}
I(\phi T)=&\det(\phi u,\phi v,\phi w)\, \textrm{Anti} (\phi (u,v,w)\textbf{e} )\\
=& \det\phi \,\det(u,v,w) \,\phi^{-t} \,\textrm{Anti}((u,v,w)\textbf{e}) \,\phi^{-1}\\
=& \phi^{-t} \,I(T)\,\phi^{-1}.
\end{align*}
Let $T \in \mathcal{T}^3$ with dim $T\leq 2$, we have $I(\phi T)=\phi^{-t} \,I(T)\,\phi^{-1}=\textbf{0}$ for all $\phi \in \textrm{SL}(3)$.
\qed

\begin{lem}\label{lem2.2}
The function $I$ is a simple $\textrm{SL}(3)$ contravariant valuation on $\mathcal{P}_o^3$.
\end{lem}

\emph{Proof.} \
Since the function $I:\mathcal{T}^3 \rightarrow \mathbb{M}^3$ can be extended to a valuation on finite unions of simplices in $\mathbb{R}^3$ that have one vertex at the origin. For $T_1, \ldots, T_m \in \mathcal{T}^3$, it follows from the inclusion-exclusion principle that
\begin{align*}
I(T_1\cup \cdots\cup T_m)=\sum_{j=1}^{m} (-1)^{j-1} \sum_{1\leq i_1<\cdots <i_j \leq m}I(T_{i_1}\cap \cdots\cap T_{i_j}).
\end{align*}
Triangulate $P \in\mathcal{P}_o^3$ into simplices $T_1, \ldots, T_l$ with one vertex at the origin. Therefore, by the simplicity of $I$ on $\mathcal{T}^3$ we derive
\begin{align*}
I(P)=I(T_1\cup \cdots\cup T_l)=\sum_{i=1}^{l} I(T_{i}).
\end{align*}
Combined with Lemma \ref{lem2.1}, it follows that $I$ is a simple $\textrm{SL}(3)$ contravariant valuation on $\mathcal{P}_o^3$.
\qed

\par  For $\nu=(1-\lambda)e_1-\lambda e_2$ and $\lambda \in (0,1)$, let $H$ be the hyperplane through the origin with normal
vector $\nu$. Write the two half-spaces bounded by $H$ in a such way of $H^+$ and $H^-$:
\begin{align*}
H^+=\{x\in \mathbb{R}^n: x\cdot ((1-\lambda)e_1-\lambda e_2)\geq 0\}
\end{align*}
and
\begin{align*}
H^-=\{x\in \mathbb{R}^n: x\cdot ((1-\lambda)e_1-\lambda e_2)\leq 0\}.
\end{align*}
This hyperplane induces the series of triangulations of $rT^n$ as well as $r\tilde{T}^n$ due to the following definition.

\begin{defn}\label{defn2}
For $\lambda \in (0,1)$, define the linear transform $\phi_1 \in \textrm{GL}(n)$ by
\begin{equation*}
\phi_1 e_1=\lambda e_1+(1-\lambda)e_2,~~~\phi_1 e_2=e_2,~~~\phi_1 e_j=e_j,\quad where~j=3,\ldots,n,
\end{equation*}
and $\phi_2 \in \textrm{GL}(n)$ by
\begin{equation*}
\phi_2 e_1=e_1,~~~\phi_2 e_2=\lambda e_1+(1-\lambda)e_2,~~~\phi_2 e_j=e_j,\quad where~j=3,\ldots,n.
\end{equation*}
\end{defn}
\par It is clear that
\begin{equation}\label{2.1}
sT^n\cap H^+=\phi_2 sT^n,~~~sT^n\cap H^-=\phi_1 sT^n,~~~and~~~~~~ sT^n\cap H=\phi_1 sT^{n-1},
\end{equation}
for every $s>0$. Similarly,
\begin{equation}\label{2.2}
s\tilde{T}^n\cap H^+=\phi_2 s\tilde{T}^n,~~~s\tilde{T}^n\cap H^-=\phi_1 s\tilde{T}^n,~~~and~~~~~~ s\tilde{T}^n\cap H=\phi_1 s\tilde{T}^{n-1},
\end{equation}
for every $s>0$.

\par The well-known solution of Cauchy's functional equation is one of the main ingredients in our proofs. Since we do not assume continuity, also functionals that depend on solutions $f: [0,\infty)\rightarrow \mathbb{R}$ of Cauchy's functional equation
\begin{align*}
f(x+y)=f(x)+f(y)
\end{align*}
for all $x,y \in [0,\infty)$. If $f$ is a measurable function, then $f$ is  linear.

\section{$\textrm{SL}(n)$ Contravariant Valuations on $\mathcal{P}_o^n$}

\subsection{The two-dimensional case.}
\begin{lem}\label{lem3.1.1}
If $\mu: \mathcal{Q}^2 \rightarrow \mathbb{M}^2$ is an $\textrm{SL}(2)$ equivariant valuation. Then
\begin{equation*}
\omega(P)=\psi_{\pi/ 2} \,\mu(P)\, \psi_{\pi/ 2}^{-1}
\end{equation*}
is an $\textrm{SL}(2)$ contravariant valuation, where $\psi_{\pi/ 2}$ is described by the rotation by the angle $\pi /2$ as
$
\begin{pmatrix}
0 & -1 \\
1 & 0
\end{pmatrix}.
$
\end{lem}
\emph{Proof.} \
Since $\mu$ is equivariant, then for every $\phi \in SL(2)$ we have
\begin{align}\label{3.1}
\omega(\phi P)=&\psi_{\pi/ 2} \,\mu(\phi P)\, \psi_{\pi/ 2}^{-1}\notag\\
=&\psi_{\pi/ 2} \,\phi\, \mu( P)\, \phi^t\, \psi_{\pi/ 2}^{-1}\notag\\
=&\psi_{\pi/ 2} \,\phi \,\psi_{\pi/ 2}^{-1} \,\psi_{\pi/ 2} \,\mu( P)\, \psi_{\pi/ 2}^{-1} \,\psi_{\pi/ 2} \, \phi^t\, \psi_{\pi/ 2}^{-1}\notag\\
=&\psi_{\pi/ 2}\, \phi\, \psi_{\pi/ 2}^{-1} \,\omega(P) \,\psi_{\pi/ 2}\, \phi^t\, \psi_{\pi/ 2}^{-1}
\end{align}
Let
$
\phi=\begin{pmatrix}
a & b \\
c & d
\end{pmatrix}
$. Therefore,
\begin{align*}
\psi_{\pi/ 2}\,\phi\,\psi_{\pi/ 2}^{-1}=\begin{pmatrix}
0 & -1 \\
1 & 0
\end{pmatrix}\begin{pmatrix}
a & b \\
c & d
\end{pmatrix}\begin{pmatrix}
0 & 1 \\
-1 & 0
\end{pmatrix}=\begin{pmatrix}
d & -c \\
-b & a
\end{pmatrix}=\phi^{-t}
\end{align*}
and
\begin{align*}
\psi_{\pi/ 2}\,\phi^t \,\psi_{\pi/ 2}^{-1}=\begin{pmatrix}
0 & -1 \\
1 & 0
\end{pmatrix}\begin{pmatrix}
a & c \\
b & d
\end{pmatrix}\begin{pmatrix}
0 & 1 \\
-1 & 0
\end{pmatrix}=\begin{pmatrix}
d & -b \\
-c & a
\end{pmatrix}=\phi^{-1}.
\end{align*}
Thus, (\ref{3.1}) is equivalent to
\begin{align*}
\omega(\phi P)=\phi^{-t}\,\omega(P)\,\phi^{-1}
\end{align*}
for every $\phi \in \textrm{SL}(2)$.
\qed


\begin{lem}\label{lem3.1.2}(\cite{Wang})
A function $\mu: \mathcal{P}_o^2 \rightarrow \mathbb{M}^2$ is an $\textrm{SL}(2)$ equivariant valuation if and only if there exist constants $c_1, c_2, c_3, c_4 \in \mathbb{R}$ and solutions of Cauchy's functional equation $\alpha, \beta: [0, \infty)\rightarrow \mathbb{R}$ such that
\begin{align}
\mu(P)=c_1\,M(P)+c_2\,E(P)+c_3\,F(P)+H_\alpha(P)+G_\beta(P)+c_4\,\rho_{\pi/ 2}
\end{align}
for every $P\in \mathcal{P}_o^2 $, where $\rho_{\pi/ 2}$ is the counterclockwise rotation in $\mathbb{R}^2$ by the angle $\pi/ 2$.
\end{lem}

\par Now, combining Lemma \ref{lem3.1.1} and Lemma \ref{lem3.1.2}, we obtain  Theorem \ref{theo3}.

\subsection{The three-dimensional case.}

\begin{lem}\label{lem3.2.1}
If $\mu: \mathcal{P}_o^3 \rightarrow \mathbb{M}^3$ is an $\textrm{SL}(3)$ contravariant valuation, then $\mu(\{o\})=\textbf{0}$.
\end{lem}

\emph{Proof.} \
Set
$\mu(\{o\})=
\begin{pmatrix}
a_{11} & a_{12} & a_{13} \\
a_{21} & a_{22} & a_{23} \\
a_{31} & a_{22} & a_{33}
\end{pmatrix}
$ and
$
\psi_1=\begin{pmatrix}
s^2 & 0 & 0 \\
0 & \frac{1}{s} & 0 \\
0 & 0 & \frac{1}{s}
\end{pmatrix}
$ for $s\neq 0$. The $\textrm{SL}(3)$ contravariance of $\mu$ implies
\begin{align*}
\mu(\{o\})=\mu(\psi_1\{o\})= &\psi_1^{-t} \mu(\{o\}) \psi_1^{-1}\\
=& \begin{pmatrix}
\frac{1}{s^2} & 0 & 0 \\
0 & s & 0 \\
0 & 0 & s
\end{pmatrix}\begin{pmatrix}
a_{11} & a_{12} & a_{13} \\
a_{21} & a_{22} & a_{23} \\
a_{31} & a_{22} & a_{33}
\end{pmatrix}\begin{pmatrix}
\frac{1}{s^2} & 0 & 0 \\
0 & s & 0 \\
0 & 0 & s
\end{pmatrix}\\
=&\begin{pmatrix}
a_{11}/s^4 & a_{12}/s & a_{13}/s \\
a_{21}/s & s^2 a_{22} & s^2a_{23} \\
a_{31}/s & s^2a_{22} & s^2a_{33}
\end{pmatrix}\\
=&\begin{pmatrix}
a_{11} & a_{12} & a_{13} \\
a_{21} & a_{22} & a_{23} \\
a_{31} & a_{22} & a_{33}
\end{pmatrix}.
\end{align*}
In the last step, we apply the arbitrariness of $s$ to obtain $a_{ij}=0$, where $1\leq i,j \leq 3$. This yields $\mu(\{o\})=\textbf{0}$.
\qed

\begin{lem}\label{lem3.2.2}
If $\mu: \mathcal{P}_o^3 \rightarrow \mathbb{M}^3$ is an $\textrm{SL}(3)$ contravariant valuation, then there exists a constant $t_1 \in \mathbb{R}$ such that
\begin{equation*}
\mu(r T^1)= t_1 r\begin{pmatrix}
0 & 0 & 0 \\
0 & 0 & 1 \\
0 & -1 & 0
\end{pmatrix}
\end{equation*}
for $r\neq 0$.
\end{lem}

\emph{Proof.} \
Set $\mu(T^1)=\begin{pmatrix}
A_{11} & A_{12} \\
A_{21} & A_{22}
\end{pmatrix}$, and $\psi_2=\begin{pmatrix}
1 & B_{12} \\
\textbf{0} & B_{22}
\end{pmatrix}$ for any pair of matrices $B_{12}$ and $B_{22}$ with $B_{22} \in \textrm{SL}(2)$. The $\textrm{SL(3)}$ contravariance of $\mu$ implies that
\begin{align}\label{3.3}
\mu(T^1)=\mu(\psi_2T^1)=& \psi_2^{-t} \mu(T^1) \psi_2^{-1}\notag\\
=&\begin{pmatrix}
1 & \textbf{0} \\
-B_{22}^{-t} B_{12}^{t} & B_{22}^{-t}
\end{pmatrix}\begin{pmatrix}
A_{11} & A_{12} \\
A_{21} & A_{22}
\end{pmatrix}\begin{pmatrix}
1 & -B_{12} B_{22}^{-1} \\
\textbf{0} & B_{22}^{-1}
\end{pmatrix}\notag\\
=&\left(\begin{array}{cl}
A_{11} & -A_{11} B_{12} B_{22}^{-1}+A_{12} B_{22}^{-1} \\
\\
-B_{22}^{-t} B_{12}^{t} A_{11}+B_{22}^{-t} A_{21} & B_{22}^{-t} B_{12}^{t} A_{11} B_{12} B_{22}^{-1}-B_{22}^{-t} A_{21} B_{12} B_{22}^{-1} \\
& -B_{22}^{-t} B_{12}^{t} A_{12} B_{22}^{-1}+B_{22}^{-t} A_{22} B_{22}^{-1}
\end{array}\right)
\end{align}
Setting $B_{12}=\textbf{0}$ in (\ref{3.3}), we have
\begin{equation}\label{3.4}
A_{12}B_{22}^{-1}=A_{12},B_{22}^{-t} A_{21}=A_{21}~\textrm{and}~B_{22}^{-t} A_{22} B_{22}^{-1}=A_{22}
\end{equation}
for $B_{22} \in \textrm{SL}(2)$. Let $B_{22}=\begin{pmatrix}
s & 0 \\
0 & \frac{1}{s}
\end{pmatrix}$ for $s\neq 0$, then (\ref{3.4}) becomes
\begin{equation}\label{3.5}
A_{12}=A_{21}=\textbf{0}~\textrm{and}~A_{22}=\begin{pmatrix}
0 & x \\
y & 0
\end{pmatrix}
\end{equation}
for any $x, y \in \mathbb{R}$. Inserting (\ref{3.5}) into (\ref{3.3}) we have
\begin{equation}\label{3.6}
A_{11} B_{12} B_{22}^{-1}=B_{22}^{-t} B_{12}^{t} A_{11}=\textbf{0}.
\end{equation}
For $B_{12}=\begin{pmatrix}
0 & 1
\end{pmatrix}$  and $B_{22}=I_{22},$ which is  the $(2\times 2)$ identity matrix, it follows from (\ref{3.6})  that

\begin{equation}\label{3.7}
A_{11}=\textbf{0}.
\end{equation}
Combining (\ref{3.3}), (\ref{3.5}) and (\ref{3.7}), then we have
\begin{equation*}
B_{22}^{-t}A_{22} B_{22}^{-1}=A_{22}.
\end{equation*}
Let $B_{22}=\begin{pmatrix}
1 & -1 \\
0 & 1
\end{pmatrix}$, then there exists a constant $t_1 \in \mathbb{R}$ such that $A_{22}=\begin{pmatrix}
0 & t_1 \\
-t_1 & 0
\end{pmatrix}$.
Therefore, applying $\psi_{3}=\begin{pmatrix}
r & 0 & 0 \\
0 & r & 0 \\
0 & 0 & \frac{1}{r^2}
\end{pmatrix}$ for $r\neq 0$, we have
\begin{align*}
\mu(r T^1)=\mu(\psi_3T^1)=& \psi_3^{-t} \mu(T^1) \psi_3^{-1}\notag\\
=&\begin{pmatrix}
\frac{1}{r} & 0 & 0 \\
0 & \frac{1}{r} & 0 \\
0 & 0 & r^2
\end{pmatrix}\begin{pmatrix}
0 & 0 & 0 \\
0 & 0 & t_1 \\
0 & -t_1 & 0
\end{pmatrix}\begin{pmatrix}
\frac{1}{r} & 0 & 0 \\
0 & \frac{1}{r} & 0 \\
0 & 0 & r^2
\end{pmatrix}\\
=&t_1 r\begin{pmatrix}
0 & 0 & 0 \\
0 & 0 & 1 \\
0 & -1 & 0
\end{pmatrix}.
\end{align*}
\qed

\begin{lem}\label{lem3.2.3}
If $\mu: \mathcal{P}_o^3 \rightarrow \mathbb{M}^3$ is an $\textrm{SL}(3)$ contravariant valuation, then there exists a constant $t_1 \in \mathbb{R}$ such that
\begin{equation*}
\mu(r T^2)= \frac{t_1 r}{2}\begin{pmatrix}
0 & 0 & -1 \\
0 & 0 & 1 \\
1 & -1 & 0
\end{pmatrix}
\end{equation*}
for $r\neq 0$.
\end{lem}

\emph{Proof.} \
Set $\mu(T^2)=\begin{pmatrix}
A_{22} & A_{21} \\
A_{12} & A_{11}
\end{pmatrix}$ and $\psi_4=\begin{pmatrix}
I_{22} & B_{21} \\
\textbf{0} & 1
\end{pmatrix}$ for any matrix $B_{21}$. The $\textrm{SL}(3)$ contravariance of $\mu$ implies that
\begin{align}\label{3.8}
\mu(T^2)=\mu(\psi_4T^2)=& \psi_4^{-t} \mu(T^2) \psi_4^{-1}\notag\\
=&\begin{pmatrix}
I_{22} & \textbf{0} \\
-B_{21}^{t} & 1
\end{pmatrix}\begin{pmatrix}
A_{22} & A_{21} \\
A_{12} & A_{11}
\end{pmatrix}\begin{pmatrix}
I_{22} & -B_{21} \\
\textbf{0} & 1
\end{pmatrix}\notag\\
=&\begin{pmatrix}
A_{22} & -A_{22}B_{21}+A_{21} \\
-B_{21}^{t}A_{22}+A_{12} & B_{21}^{t}A_{22}B_{21}-A_{12}B_{21}-B_{21}^{t}A_{21}+A_{11}
\end{pmatrix}.
\end{align}
This yields
\begin{equation}\label{3.9}
B_{21}^{t}A_{22}=A_{22}B_{21}=\textbf{0}.
\end{equation}
For $B_{21}=\begin{pmatrix}
1 \\
0
\end{pmatrix}$, then we have $A_{22}=\textbf{0}$ and $A_{12}B_{21}+B_{21}^{t}A_{21}=0.$
 Setting  $A_{12}=\begin{pmatrix}
a_1 & a_2
\end{pmatrix},$ $A_{21}=\begin{pmatrix}
a_3 \\
a_4
\end{pmatrix}$ and $B_{21}=\begin{pmatrix}
c \\
d
\end{pmatrix}$ for every $c,d \in \mathbb{R}$, it follows from (\ref{3.9}) that
\begin{equation*}
\begin{pmatrix}
a_1 & a_2
\end{pmatrix}\begin{pmatrix}
c \\
d
\end{pmatrix}+\begin{pmatrix}
c & d
\end{pmatrix}\begin{pmatrix}
a_3 \\
a_4
\end{pmatrix}=c(a_1+a_3)+d(a_2+a_4)=0
\end{equation*}
for any $c,d \in \mathbb{R}$. It leads to $a_1+a_3=a_2+a_4=0$. Next, for
$
\psi_5=\begin{pmatrix}
0 & 1 & 0 \\
1 & 0 & 1 \\
0 & 0 & -1
\end{pmatrix}
\in \textrm{SL}(3),$ and by $\psi_5T^2=T^2$ and the $\textrm{SL}(3)$ contravariance of $\mu$, we have
\begin{align*}
\begin{pmatrix}
0 & 0 & a_2 \\
0 & 0 & a_1 \\
-a_2 & -a_1 & A_{11}
\end{pmatrix}
=\begin{pmatrix}
0 & 0 & -a_1 \\
0 & 0 & -a_2 \\
a_1 & a_2 & A_{11}
\end{pmatrix}.
\end{align*}
Thus, $a_1=-a_2$, and there exist constants $t_2, t_3 \in \mathbb{R}$ such that $\mu(T^2)=\begin{pmatrix}
0 & 0 & -t_2 \\
0 & 0 & t_2 \\
t_2 & -t_2 & t_3
\end{pmatrix}$.

On the other hand, using $\psi_3$, it follows that
\begin{align}\label{3.10}
\mu(r T^2)=\mu(\psi_3 T^2)=& \psi_3^{-t} \mu(T^2) \psi_3^{-1}\notag\\
=&\begin{pmatrix}
\frac{1}{r} & 0 & 0 \\
0 & \frac{1}{r} & 0 \\
0 & 0 & r^2
\end{pmatrix}\begin{pmatrix}
0 & 0 & -t_2 \\
0 & 0 & t_2 \\
t_2 & -t_2 & t_3
\end{pmatrix}\begin{pmatrix}
\frac{1}{r} & 0 & 0 \\
0 & \frac{1}{r} & 0 \\
0 & 0 & r^2
\end{pmatrix}\notag\\
=&\begin{pmatrix}
0 & 0 & -t_2r \\
0 & 0 & t_2r \\
t_2r & -t_2r & t_3r^4
\end{pmatrix}.
\end{align}
\par Let $\phi_1, \phi_2 \in\textrm{GL}(3)$ be defined as in Definition \ref{defn2}. Since
\begin{equation*}
T^2\cap H^+=\phi_2 T^2,~~~T^2\cap H^-=\phi_1 T^2,~~~\textrm{and}~~~~~~ T^2\cap H=\phi_1 T^1,
\end{equation*}
then the valuation property of $\mu$ gives
\begin{equation*}
\mu(T^2)+\mu(\phi_1 T^1)=\mu(\phi_1 T^2)+\mu(\phi_2 T^2).
\end{equation*}
Note that $\phi_1/ \lambda^{\frac{1}{3}}$, $\phi_2/ (1-\lambda)^{\frac{1}{3}}$ belong to $\textrm{SL}(3)$, then by the $\textrm{SL}(3)$ contravariance of $\mu$ we have
\begin{align}\label{3.11}
\mu(T^2)+ \lambda^{\frac{2}{3}} {\phi_1}^{-t}\mu\left( \lambda^{\frac{1}{3}} T^1\right) {\phi_1}^{-1}
=&\lambda^{\frac{2}{3}} {\phi_1}^{-t}\mu\left( \lambda^{\frac{1}{3}} T^2\right) {\phi_1}^{-1}\\ \nonumber
&+(1-\lambda)^{\frac{2}{3}} {\phi_2}^{-t}\mu\left( (1-\lambda)^{\frac{1}{3}} T^2\right) {\phi_2}^{-1}.
\end{align}
Combining Lemma \ref{lem3.2.2} and (\ref{3.10}), (\ref{3.11}) becomes
\begin{align*}
&\begin{pmatrix}
0 & 0 & -t_2 \\
0 & 0 & t_2 \\
t_2 & -t_2 & t_3
\end{pmatrix}+t_1\begin{pmatrix}
0 & 0 & \lambda-1 \\
0 & 0 & \lambda \\
1-\lambda & -\lambda & 0
\end{pmatrix}\\
=&\lambda^{\frac{2}{3}} \begin{pmatrix}
\frac{1}{\lambda} & \frac{\lambda-1}{\lambda} & 0 \\
0 & 1 & 0 \\
0 & 0 & 1
\end{pmatrix}\begin{pmatrix}
0 & 0 & -\lambda^{\frac{1}{3}}t_2 \\
0 & 0 & \lambda^{\frac{1}{3}}t_2 \\
\lambda^{\frac{1}{3}}t_2 & -\lambda^{\frac{1}{3}}t_2 & \lambda^{\frac{4}{3}}t_3
\end{pmatrix}\begin{pmatrix}
\frac{1}{\lambda} & 1 & 0 \\
\frac{\lambda-1}{\lambda} & 1 & 0 \\
0 & 0 & 1
\end{pmatrix}\\
&+(1-\lambda)^{\frac{2}{3}} \begin{pmatrix}
1 & 0 & 0 \\
-\frac{\lambda}{1-\lambda} & \frac{1}{1-\lambda} & 0 \\
0 & 0 & 1
\end{pmatrix}\begin{pmatrix}
0 & 0 & -(1-\lambda)^{\frac{1}{3}}t_2 \\
0 & 0 & (1-\lambda)^{\frac{1}{3}}t_2 \\
(1-\lambda)^{\frac{1}{3}}t_2 & -(1-\lambda)^{\frac{1}{3}}t_2 & (1-\lambda)^{\frac{4}{3}}t_3
\end{pmatrix}\begin{pmatrix}
1 & -\frac{\lambda}{1-\lambda} & 0 \\
0 & \frac{1}{1-\lambda} & 0 \\
0 & 0 & 1
\end{pmatrix}\\
=&\begin{pmatrix}
0 & 0 & (\lambda-2)t_2 \\
0 & 0 & \lambda t_2 \\
(2-\lambda)t_2 & -\lambda t_2 & \lambda^2 t_3
\end{pmatrix}+\begin{pmatrix}
0 & 0 & (\lambda-1)t_2 \\
0 & 0 & (\lambda+1) t_2 \\
(1-\lambda)t_2 & -(\lambda+1) t_2 & (1-\lambda)^2 t_3
\end{pmatrix}\\
=&\begin{pmatrix}
0 & 0 & (2\lambda-3)t_2 \\
0 & 0 & (2\lambda+1) t_2 \\
(3-2\lambda)t_2 & -(2\lambda+1) t_2 & (2\lambda^2-2\lambda+1) t_3
\end{pmatrix}
\end{align*}
for every $0<\lambda <1$. Therefore,
\begin{equation*}
t_2=\frac{t_1}{2}~\textrm{and}~t_3=0,
\end{equation*}
which yields
\begin{align*}
\mu(r T^2)=\frac{t_1 r}{2}\begin{pmatrix}
0 & 0 & -1 \\
0 & 0 & 1 \\
1 & -1 & 0
\end{pmatrix}.
\end{align*}
\qed

\begin{lem}\label{lem3.2.4}
If $\mu: \mathcal{P}_o^3 \rightarrow \mathbb{M}^3$ is an $\textrm{SL}(3)$ contravariant valuation, then $\mu$ is simple. Moreover, there exists a constant $c \in \mathbb{R}$ such that
\begin{align*}
\mu(r^{\frac{1}{3}} T^3)=c \,I(r^{\frac{1}{3}} T^3)+2L_\xi(r^{\frac{1}{3}} T^3)
\end{align*}
for $r> 0$.
\end{lem}

\emph{Proof.} \
By the valuation property of $\mu$, it follows that
\begin{equation*}
\mu(r^{\frac{1}{3}} T^3)+\mu(\phi_1 r^{\frac{1}{3}}T^2)=\mu(\phi_1 r^{\frac{1}{3}} T^3)+\mu(\phi_2 r^{\frac{1}{3}}T^3).
\end{equation*}
Thus
\begin{align}\label{3.12}
\mu(r^{\frac{1}{3}} T^3)+ \lambda^{\frac{2}{3}} {\phi_1}^{-t}\mu\left( (\lambda r)^{\frac{1}{3}} T^2\right) {\phi_1}^{-1}
=&\lambda^{\frac{2}{3}} {\phi_1}^{-t}\mu\left( (\lambda r)^{\frac{1}{3}} T^3\right) {\phi_1}^{-1}\\
&+(1-\lambda)^{\frac{2}{3}} {\phi_2}^{-t}\mu\left( ((1-\lambda)r)^{\frac{1}{3}} T^3\right) {\phi_2}^{-1}.\nonumber
\end{align}
Setting $\lambda=a/a+b$ and $r=a+b$ for $a, b >0$. Combining with Lemma \ref{lem3.2.3}, it follows from (\ref{3.12}) that
\begin{align}\label{3.13}
(a+b)^{\frac{2}{3}}\mu\left((a+b)^{\frac{1}{3}} T^3\right) &+ \frac{t_1}{2}\begin{pmatrix}
0 & 0 & -(a+2b) \\
0 & 0 & a \\
a+2b & -a & 0
\end{pmatrix}\\  \nonumber
=&a^{\frac{2}{3}} {\phi_1}^{-t}\mu(a^{\frac{1}{3}} T^3){\phi_1}^{-1}
+b^{\frac{2}{3}} {\phi_2}^{-t}\mu(b^{\frac{1}{3}} T^3){\phi_2}^{-1}. \nonumber
\end{align}
For $x^{\frac{2}{3}} \mu(x^{\frac{1}{3}} T^3)=\begin{pmatrix}
g_{11}(x) & g_{12}(x) & g_{13}(x) \\
g_{21}(x) & g_{22}(x) & g_{23}(x) \\
g_{31}(x) & g_{32}(x) & g_{33}(x)
\end{pmatrix}$ and  $\psi_6=\begin{pmatrix}
0 & 0 & 1 \\
1 & 0 & 0 \\
0 & 1 & 0
\end{pmatrix}$, we have
\begin{align*}
\mu(x^{\frac{1}{3}} T^3)=\mu(\psi_6 x^{\frac{1}{3}} T^3)=& \psi_6^{-t} \mu(x^{\frac{1}{3}} T^3) \psi_6^{-1}\\
=&\begin{pmatrix}
0 & 0 & 1 \\
1 & 0 & 0 \\
0 & 1 & 0
\end{pmatrix}x^{-\frac{2}{3}}\begin{pmatrix}
g_{11}(x) & g_{12}(x) & g_{13}(x) \\
g_{21}(x) & g_{22}(x) & g_{23}(x) \\
g_{31}(x) & g_{32}(x) & g_{33}(x)
\end{pmatrix}\begin{pmatrix}
0 & 0 & 1 \\
1 & 0 & 0 \\
0 & 1 & 0
\end{pmatrix}\\
=&x^{-\frac{2}{3}}\begin{pmatrix}
g_{33}(x) & g_{31}(x) & g_{32}(x) \\
g_{13}(x) & g_{11}(x) & g_{12}(x) \\
g_{23}(x) & g_{21}(x) & g_{22}(x)
\end{pmatrix}\\
=&x^{-\frac{2}{3}}\begin{pmatrix}
g_{11}(x) & g_{12}(x) & g_{13}(x) \\
g_{21}(x) & g_{22}(x) & g_{23}(x) \\
g_{31}(x) & g_{32}(x) & g_{33}(x)
\end{pmatrix}.
\end{align*}
It leads to
\begin{equation*}
g_{11}=g_{22}=g_{33}, g_{12}=g_{23}=g_{31}~\textrm{and}~g_{13}=g_{21}=g_{32}.
\end{equation*}
Therefore, (\ref{3.13}) can be rewritten as
\begin{align*}
&\begin{pmatrix}
g_{1}(a+b) & g_{2}(a+b) & g_{3}(a+b) \\
g_{3}(a+b) & g_{1}(a+b) & g_{2}(a+b) \\
g_{2}(a+b) & g_{3}(a+b) & g_{1}(a+b)
\end{pmatrix}+\frac{t_1}{2}\begin{pmatrix}
0 & 0 & -(a+2b) \\
0 & 0 & a \\
a+2b & -a & 0
\end{pmatrix}\\
=&\begin{pmatrix}
\frac{(a+b)^{2}+b^2}{a^{2}} g_{1}(a)-\frac{(a+b) b}{a^{2}} g_{3}(a)-\frac{(a+b) b}{a^{2}} g_{2}(a) & \frac{a+b}{a} g_{2}(a)-\frac{b}{a} g_{1}(a) & \frac{a+b}{a} g_{3}(a)-\frac{b}{a} g_{2}(a) \\
\frac{a+b}{a} g_{3}(a)-\frac{b}{a} g_{1}(a) & g_{1}(a) & g_{2}(a) \\
\frac{a+b}{a} g_{2}(a)-\frac{b}{a} g_{3}(a) & g_{3}(a) & g_{1}(a)
\end{pmatrix}\\
&+\begin{pmatrix}
g_{1}(b) & \frac{a+b}{b} g_{2}(b)-\frac{a}{b} g_{1}(b) & g_{3}(b) \\
\frac{a+b}{b} g_{3}(b)-\frac{a}{b} g_{1}(b) & \frac{(a+b)^{2}+a^{2}}{b^{2}} g_{1}(b)-\frac{(a+b) a}{b^{2}} g_{3(b)}-\frac{(a+b) a}{b^{2}} g_{2}(b) & \frac{a+b}{b} g_{2}(b)-\frac{a}{b} g_{3}(b)  \\
g_{2}(b) & \frac{a+b}{b} g_{3}(b)-\frac{a}{b} g_{2}(b) & g_{1}(b)
\end{pmatrix}.
\end{align*}
Thus,
\begin{equation}\label{3.14}
g_{1}(a+b)=\frac{(a+b)^{2}+b^2}{a^{2}} g_{1}(a)-\frac{(a+b) b}{a^{2}} g_{3}(a)-\frac{(a+b) b}{a^{2}} g_{2}(a)+g_{1}(b),
\end{equation}
\begin{equation}\label{3.15}
g_{3}(a+b)=\frac{a+b}{a} g_{3}(a)-\frac{b}{a} g_{1}(a)+\frac{a+b}{b} g_{3}(b)-\frac{a}{b} g_{1}(b),
\end{equation}
\begin{equation}\label{3.16}
g_{2}(a+b)+\frac{t_1(a+2b)}{2}=\frac{a+b}{a} g_{2}(a)-\frac{b}{a} g_{3}(a)+g_{2}(b),
\end{equation}
\begin{equation}\label{3.17}
g_{2}(a+b)=\frac{a+b}{a} g_{2}(a)-\frac{b}{a} g_{1}(a)+\frac{a+b}{b} g_{2}(b)-\frac{a}{b} g_{1}(b),
\end{equation}
\begin{equation}\label{3.18}
g_{1}(a+b)=g_{1}(a)+\frac{(a+b)^{2}+a^{2}}{b^{2}} g_{1}(b)-\frac{(a+b) a}{b^{2}} g_{3(b)}-\frac{(a+b) a}{b^{2}} g_{2}(b),
\end{equation}
\begin{equation}\label{3.19}
g_{3}(a+b)-\frac{t_1 a}{2}=g_{3}(a)+\frac{a+b}{b} g_{3}(b)-\frac{a}{b} g_{2}(b),
\end{equation}
\begin{equation}\label{3.20}
g_{3}(a+b)-\frac{t_1(a+2b)}{2}=\frac{a+b}{a} g_{3}(a)-\frac{b}{a} g_{2}(a)+g_{3}(b),
\end{equation}
\begin{equation}\label{3.21}
g_{2}(a+b)+\frac{t_1 a}{2}=g_{2}(a)+\frac{a+b}{b} g_{2}(b)-\frac{a}{b} g_{3}(b),
\end{equation}
\begin{equation}\label{3.22}
g_{1}(a+b)=g_{1}(a)+ g_{1}(b).
\end{equation}
By (\ref{3.15}) and (\ref{3.17}), we obtain
\begin{equation*}
g_{2}(a+b)-g_{3}(a+b)=\frac{a+b}{a}\left(g_{2}(a)-g_{3}(a)\right)+\frac{a+b}{b}\left(g_{2}(b)-g_{3}(b)\right),
\end{equation*}
which is equivalent to
\begin{equation}\label{3.23}
g_{2}(x)-g_{3}(x)=\eta(x)x,
\end{equation}
where $\eta: [0, \infty) \rightarrow \mathbb{R}$ is a solution of Cauchy's functional equation. Applying (\ref{3.16}), (\ref{3.21}) and (\ref{3.23}), we have
\begin{equation}\label{3.24}
t_1=\eta(a)-\frac{a}{b}\eta(b).
\end{equation}
It follows from (\ref{3.20}), (\ref{3.21}) and (\ref{3.23}) that
$\eta(a+b)+t_1=\eta(a)+\eta(b)$. Hence,
\begin{equation}\label{3.25}
t_1=0.
\end{equation}
Inserting (\ref{3.25}) into (\ref{3.24}), we have $\eta(a)/a=\eta(b)/b$. So there exists a constant $t_4 \in \mathbb{R}$ such that $\eta(x)=t_4 x$. Further, combining (\ref{3.16}), (\ref{3.23}) and (\ref{3.25}), we have
$g_{2}(a+b)=g_{2}(a)+g_{2}(b)+ab\,t_4$. This yields
\begin{equation}\label{3.26}
g_{2}(x)=\xi(x)+\frac{t_4 x^2}{2},
\end{equation}
where $\xi: [0, \infty) \rightarrow \mathbb{R}$ is a solution of Cauchy's functional equation. So
\begin{equation}\label{3.27}
g_{3}(x)=\xi(x)-\frac{t_4 x^2}{2}.
\end{equation}
By (\ref{3.14}), (\ref{3.26}) and (\ref{3.27}), we obtain
\begin{equation}\label{3.28}
g_{1}(a+b)=\frac{a^2+2b^2+2ab}{a^2}g_1(a)+g_1(b)-\frac{2(a+b)b}{a^2}\xi(a).
\end{equation}
Similarly, it follows from (\ref{3.18}), (\ref{3.26}) and (\ref{3.27}) that
\begin{equation}\label{3.29}
g_{1}(a+b)=g_1(a)+\frac{2a^2+b^2+2ab}{b^2}g_1(b)-\frac{2(a+b)a}{b^2}\xi(b).
\end{equation}
Combining (\ref{3.28}) and (\ref{3.29}), then we have
\begin{equation*}
\frac{g_1(a)-\xi(a)}{a^3}=\frac{g_1(b)-\xi(b)}{b^3}.
\end{equation*}
Hence, there exists a constant $h\in \mathbb{R}$ such that
\begin{equation}\label{3.30}
g_1(x)=\xi(x)+hx^3.
\end{equation}
Inserting (\ref{3.22}) into (\ref{3.14}) it concludes that $2g_1(a)-g_2(a)-g_3(a)=0$. Applying (\ref{3.26}), (\ref{3.27}) and (\ref{3.30}) to obtain $h=0$. Thus,
\begin{equation}\label{3.31}
g_1(x)=\xi(x).
\end{equation}
\par From (\ref{3.25}), Lemma \ref{lem3.2.2} and Lemma \ref{lem3.2.3}, we have  $\mu(rT^1)=\mu(rT^2)=\textbf{0}$ for every $r>0$. Since every $k$-dimensional simplex $T \in \mathcal{P}_o^3$ is an $\textrm{SL}(3)$ image of $rT^k$ for $k=1,2$, we obtain $\mu(T)=\textbf{0}$ for every simplex $T$ with $\dim T\leq 2$. Moreover, applying the triangulation of $k$-dimensional polytope $P$ into simplices as a set
of $k$-dimensional simplices $\{T_1,\ldots,T_m\}$, which have one
vertex at the origin and pairwise disjoint interiors. Then it follows from the inclusion-exclusion principle that $\mu(P)=\textbf{0}$ for every $P\in \mathcal{P}_o^3$ on $\dim T\leq 2$, that is, $\mu$ is simple.

Applying (\ref{3.26}), (\ref{3.27}) and (\ref{3.31}), it yields
\begin{align*}
\mu(r^{\frac{1}{3}} T^3)=&r^{-\frac{2}{3}}\begin{pmatrix}
\xi(r) & \xi(r)+\frac{t_4 r^2}{2} & \xi(r)-\frac{t_4 r^2}{2} \\
\xi(r)-\frac{t_4 r^2}{2} & \xi(r) & \xi(r)+\frac{t_4 r^2}{2} \\
\xi(r)+\frac{t_4 r^2}{2} & \xi(r)-\frac{t_4 r^2}{2} & \xi(r)
\end{pmatrix}\\
=&t_4 r^{\frac{4}{3}}\begin{pmatrix}
0 & \frac{1}{2} & -\frac{1}{2} \\
-\frac{1}{2} & 0 & \frac{1}{2} \\
\frac{1}{2} & -\frac{1}{2} & 0
\end{pmatrix}+ \xi (r) r^{-\frac{2}{3}}\begin{pmatrix}
1 & 1 & 1 \\
1 & 1 & 1 \\
1 & 1 & 1
\end{pmatrix}\\
=&-\frac{t_4}{2}  I(r^{\frac{1}{3}} T^3)+2L_\xi(r^{\frac{1}{3}} T^3).
\end{align*}
Set $c=-\frac{t_4}{2}$ to complete the proof.
\qed

\par Next, we give the proof of Theorem \ref{theo2}.

\textbf{Proof of Theorem \ref{theo2}.} \
By Lemma \ref{lem2.2} we obtain that (\ref{3.2.1}) determines an $\textrm{SL}(3)$ contravariant valuation on $ \mathcal{P}_o^3$. On the other hand, combining Lemma \ref{lem 2}, Lemma \ref{lem3.2.1} and Lemma \ref{lem3.2.4} we complete the reverse statement.
\qed

\subsection{The higher-dimensional case.}
First, we use the tensor representation for $A=(a_{ij}) \in \mathbb{M}$, that is
\begin{equation*}
A=\sum_{1\leq i \leq j \leq n } a_{ij} e_i\otimes e_j,
\end{equation*}
and write $a_{ij}=A(e_i, e_j)$. For every $\phi \in \textrm{GL}(n)$ and $y_1, y_2 \in \mathbb{R}^n$, we define
\begin{equation*}
(\phi^{-t} \cdot A)(y_1, y_2)=A(\phi^{-1}y_1, \phi^{-1}y_2),
\end{equation*}
which coincides with the action $\phi^{-t} A \phi^{-1}$ in Ludwig \cite{Ludwig1,Ludwig2,Ludwig3,Ludwig4} in the following way
\begin{equation*}
\phi^{-t}  \cdot A=\sum_{1\leq i \leq j \leq n } a_{ij}(\phi^{-1} e_i) \otimes (\phi^{-1} e_j)=\sum_{1\leq i \leq j \leq n } a_{ij} \phi^{-t} (e_i\otimes e_j) \phi^{-1}=\phi^{-t} A \phi^{-1}.
\end{equation*}

\begin{lem}\label{lem3.3.1}\cite{Ma-Wang}
Let $n\geq 4$. If $\mu: \mathcal{P}_o^n \rightarrow \mathbb{M}_e^n$ is an $\textrm{SL}(n)$ contravariant valuation, then $\mu$ is simple.
\end{lem}


\begin{lem}\label{lem3.3.2}
Let $n\geq 4$. If $\mu: \mathcal{P}_o^n \rightarrow \mathbb{M}^n$ is an $\textrm{SL}(n)$ contravariant valuation, then $\mu$ is simple.
\end{lem}

\emph{Proof.} \
Due to Lemma \ref{lem 2} and Lemma \ref{lem3.3.1}, it suffices to prove that $\mu(r T^k)$ are symmetric for all $0\leq k \leq n-1$ and $r>0$.
\par  First, we show that $\mu(\{o\})$ is symmetric. Using the row interchanges of the $i$-th row and the $j$-th row of $I_n$,  as well as the $k$-th row and the $l$-th row,  we have the following transformation
\begin{align*}
\begin{array}{cl}
\qquad \ \ \ (i) \quad \ \, (j) \quad \ \ (k) \quad  \ \, (l) \qquad \qquad \
\end{array}\\
\tau_0=
\setlength{\arraycolsep}{2.5pt}
\renewcommand{\arraystretch}{0.4}
\addtocounter{MaxMatrixCols}{10}
\begin{pmatrix}
1  \\
  & \ddots  \\
  & & 0 & \cdots & 1 \\
  & & \vdots & \ddots & \vdots  \\
  & & 1 & \cdots & 0 \\
  & & & & &\ddots  \\
  & & & & & & 0 & \cdots & 1 \\
  & & & & & & \vdots & \ddots & \vdots  \\
  & & & & & & 1 & \cdots & 0 \\
  & & & & & & & & & \ddots \\
  & & & & & & & & & & 1 \\
\end{pmatrix}
\begin{array}{cl}
(i) \\
   \\
   \\
(j)\\
 \\
 \\
 \\
(k) \\
  \\
  \\
(l)\\
\end{array}
\end{align*}
for $1\leq i,j,k,l\leq n$. Note that $\tau_0\in \textrm{SL}(n)$ and $\tau_0^{-1}=\tau_0$, then we have
\begin{equation*}
\mu(\{o\})(e_i,e_j)=\mu(\tau_0\{o\})(e_i,e_j)=\mu(\{o\})(\tau_0^{-1}e_i,\tau_0^{-1}e_j)=\mu(\{o\})(e_j,e_i)
\end{equation*}
for pairwise different numbers $1 \leq i,j\leq n$. Thus, $\mu(\{o\})$ is symmetric.

\par Second, we  consider that $\mu(r T^1)$ is symmetric. For $1\leq i <j\leq n$, multiplying both the $j$-th row and the $l$-th row of $I_n$ by $-1$, we have $\tau_1 \in \textrm{SL}(n)$ as follows:
\begin{align*}
\begin{array}{cl}
\qquad \qquad \qquad (j)  \qquad \ \ \  (l) \qquad \qquad \qquad
\end{array}\\
\tau_1=\begin{pmatrix}
1 &  &  &  &  &  &\\
  & \ddots &  &  &  &  &   \\
  &  & -1 &  &  &  & \\
  &  &  & \ddots &  &  & \\
  &  &  &  & -1 &  & \\
  &  &  &  & & \ddots & \\
  &  &  &  & &  & 1\\
\end{pmatrix}
\begin{array}{cl}
(j) \\
   \\
(l) \\
\end{array},
\end{align*}
where $1 <l\leq n$ and $l\neq i$. From $\tau_1^{-1}=\tau_1$ and the $\textrm{SL}(n)$ contravariance of $\mu$ we obtain
\begin{align*}
\mu(r T^1)(e_i,e_j)=\mu(\tau_1 r T^1)(e_i,e_j)=\mu(r T^1)(\tau_1^{-1}e_i,\tau_1^{-1}e_j)=-\mu(r T^1)(e_i,e_j),
\end{align*}
which yields $\mu(r T^1)(e_i,e_j)=0$. Similarly, it follows that $\mu(r T^1)(e_j,e_i)=0$.

\par Third, we will prove that $\mu(r T^k)$ is symmetric for $2\leq k \leq n-2$ in  three cases.
\par (i) Let $1 \leq i<j\leq k$. For $I_n$, we first use the row interchange of the $i$-th row and the $j$-th row, and then multiply the $n$-th row by $-1$. This leads to
\begin{align*}
\begin{array}{cl}
\qquad \qquad (i) \qquad   (j) \qquad \qquad \qquad \ \ \
\end{array}\\
\tau_2=
\setlength{\arraycolsep}{2.5pt}
\renewcommand{\arraystretch}{0.4}
\addtocounter{MaxMatrixCols}{10}
\begin{pmatrix}
1 &  &  &  &  &  &  &   &  &\\
  & \ddots &  &  &  &  &  &  &  &\\
  &  & 1 &  &  &  &  &  &  &\\
  &  &  & 0 & \cdots & 1 &  &  & &\\
  &  &  & \vdots & \ddots & \vdots &  &  & &\\
  &  &  & 1 & \cdots & 0 &  &  &  &\\
  &  &  &  &  &  & 1 &  &  &\\
  &  &  &  &  &  &  & \ddots &  &\\
    &  &  &  &  &  &  &  & 1 &\\
  &  &  &  &  &  &  &  &  & -1\\
\end{pmatrix}
\begin{array}{cl}
(i) \\
   \\
   \\
(j) \\
 \\
\end{array}
\end{align*}
and $\tau_2 \in \textrm{SL}(n)$. Since $\tau_2^{-1}=\tau_2$, the $\textrm{SL}(n)$ contravariance of $\mu$ gives that
\begin{align*}
\mu(r T^k)(e_i,e_j)=\mu(\tau_2 r T^k)(e_i,e_j)=\mu(r T^k)(\tau_2^{-1}e_i,\tau_2^{-1}e_j)=&\mu(r T^k)(e_j,e_i).
\end{align*}

\par (ii) Let $1 \leq i\leq k<j\leq n$. Consider the linear transform
\begin{align*}
&\begin{matrix}
\ \  \  \ \  \  \  \  \ \  \ \ \ \ (j) & \ \ \ \ & (l)
\end{matrix}\\
\tau_3=&\begin{pmatrix}
1 &  &  &  & &  & \\
  & \ddots &  &  & &  & \\
  &  & -1 &  & &  & \\
  &  &  & \ddots & &  & \\
  &  &  &  & -1 &  & \\
  &  &  &  & & \ddots & \\
  &  &  &  & &  & 1\\
\end{pmatrix}
\begin{array}{cl}
(j) \\
   \\
(l)\\
\end{array}
\end{align*}
for $k <l \leq n$. 
Since $\tau_3 \in \textrm{SL}(n)$ and $\tau_3^{-1}=\tau_3$, and by the $\textrm{SL}(n)$ contravariance of $\mu$ we obtain
\begin{align*}
\mu(r T^k)(e_i,e_j)=\mu(\tau_3 r T^k)(e_i,e_j)=\mu(r T^k)(\tau_3^{-1}e_i,\tau_3^{-1}e_j)=-\mu(r T^k)(e_i,e_j).
\end{align*}
Thus $\mu(r T^k)(e_i,e_j)=0$. Similarly, we have $\mu(r T^k)(e_j,e_i)=0$.

\par (iii) Let $k < i<j\leq n$ and
\begin{align*}
&\begin{matrix}
\ \  \  \ \  \  \  \  \ \  \ \ \ \ \ \ \ \ \ \ \ \ \ (i)
\end{matrix}\\
\tau_4=&\begin{pmatrix}
0 & 1 &  &  & &  \\
1 & 0 &  &  & &  \\
  &  & 1 &  & &  \\
  &  & &\ddots &  & &  \\
  &  &  & &-1 &  & \\
  &  &  &  & &\ddots &  \\
  &  &  &  & & &1
\end{pmatrix}
\begin{array}{cl}
   \\
   \\
(i) \\
\end{array}\in \textrm{SL}(n).
\end{align*}
By the $\textrm{SL}(n)$ contravariance of $\mu$, we have
\begin{align*}
\mu(r T^k)(e_i,e_j)=\mu(\tau_4 r T^k)(e_i,e_j)=\mu(r T^k)(\tau_4^{-1}e_i,\tau_4^{-1}e_j)=-\mu(r T^k)(e_i,e_j).
\end{align*}
It concludes that $\mu(r T^k)(e_i,e_j)=0$. Similarly, $\mu(r T^k)(e_j,e_i)=0$.

\par Next, we derive that $\mu(r T^{n-1})$ is symmetric in the following two cases.
\par (i) Let $1 \leq i<j\leq n-1$ and
\begin{align*}
\begin{array}{cl}
\quad (i) \quad \ \ (j) \qquad \qquad \qquad \qquad \qquad\ \ \
\end{array}\\
\tau_5=
\setlength{\arraycolsep}{2.5pt}
\renewcommand{\arraystretch}{0.4}
\addtocounter{MaxMatrixCols}{10}
\begin{pmatrix}
1 &  &  &  &  &  &  &   &  &\\
  & \ddots &  &  &  &  &  &  &  &\\
  &  & 1 &  &  &  &  &  &  &\\
  &  &  & 0 & \cdots & 1 &  &  & &\\
  &  &  & \vdots & \ddots & \vdots &  &  & &\\
  &  &  & 1 & \cdots & 0 &  &  &  &\\
  &  &  &  &  &  & 1 &  &  &\\
  &  &  &  &  &  &  & \ddots &  &\\
    &  &  &  &  &  &  &  & 1 &\\
  &  &  &  &  &  &  &  &  & -1\\
\end{pmatrix}
\begin{array}{cl}
(i) \\
   \\
   \\
(j) \\
 \\
\end{array}\in \textrm{SL}(n).
\end{align*}
Due to the $\textrm{SL}(n)$ contravariance of $\mu$ and $\tau_5^{-1}=\tau_5$, we obtain
\begin{align*}
\mu(r T^{n-1})(e_i,e_j)=\mu(\tau_5 r T^{n-1})(e_i,e_j)=\mu(r T^{n-1})(\tau_5^{-1}e_i,\tau_5^{-1}e_j)=&\mu(r T^{n-1})(e_j,e_i).
\end{align*}

\par (ii) Let $1 \leq i\leq n-1$. Next, we consider $\tau_6 \in \textrm{SL}(n)$ as follows:
\begin{align*}
\begin{array}{cl}
\qquad \qquad (l) \qquad   (k) \qquad \qquad \qquad \ \ \
\end{array}\\
\tau_6=
\setlength{\arraycolsep}{2.5pt}
\renewcommand{\arraystretch}{0.4}
\addtocounter{MaxMatrixCols}{10}
\begin{pmatrix}
1 &  &  &  &  &  &  &   &  &\\
  & \ddots &  &  &  &  &  &  &  &\\
  &  & 1 &  &  &  &  &  &  &\\
  &  &  & 0 & \cdots & 1 &  &  & &\\
  &  &  & \vdots & \ddots & \vdots &  &  & &\\
  &  &  & 1 & \cdots & 0 &  &  &  &\\
  &  &  &  &  &  & 1 &  &  &\\
  &  &  &  &  &  &  & \ddots &  &\\
    &  &  &  &  &  &  &  & 1 &\\
  &  &  &  &  &  &  &  &  & -1\\
\end{pmatrix}
\begin{array}{cl}
(l) \\
   \\
   \\
(k) \\
 \\
\end{array},
\end{align*}
where $1\leq l<k\leq n-1$ and $l,k\neq i$. Applying the $\textrm{SL}(n)$ contravariance of $\mu$ we have
\begin{align*}
\mu(r T^{n-1})(e_i,e_n)=\mu(\tau_6 r T^{n-1})(e_i,e_n)=\mu(r T^{n-1})(\tau_6^{-1}e_i,\tau_6^{-1}e_n)=-\mu(r T^{n-1})(e_i,e_n).
\end{align*}
This yields $\mu(r T^{n-1})(e_i,e_n)=0$. Similarly, we know $\mu(r T^{n-1})(e_n,e_i)=0$.
\qed


\textbf{Proof of Theorem \ref{theo1}.} \
Applying $\tau_0 \in \textrm{SL}(n)$  in the proof of Lemma \ref{lem3.3.2}. Since $\tau_0 r T^n=r T^n$, the $\textrm{SL}(n)$ contravariance of $\mu$ implies
\begin{equation*}
\mu(r T^n)(e_i,e_j)=\mu(\tau_0 r T^n)(e_i,e_j)=\mu(r T^n)(\tau_0^{-1}e_i,\tau_0^{-1}e_j)=\mu(r T^n)(e_j,e_i),
\end{equation*}
where $i \neq j$ and $r >0$. Hence, $\mu(r T^n)$ is symmetric. Moreover, by Lemma \ref{lem3.3.2} we obtain $\mu (r \tilde{T}^k)=\textbf{0}$ for $1\leq k \leq n-1$. Then the symmetry consumption of the images of $\mu$ can be removed in the proof of Ma and Wang \cite{Ma-Wang}. Combining with Lemma \ref{lemLYZ} we complete the proof.
\qed

\section{$\textrm{SL}(n)$ Contravariant Valuations on $\mathcal{P}^n$}

\subsection{The two-dimensional case.}

\begin{lem}\cite{Wang}\label{lem4.1.1}
A function $\mu: \mathcal{P}^2 \rightarrow \mathbb{M}^2$ is an $\textrm{SL}(2)$ equivariant valuation if and only if there exist constants $c_1, c_2, c_3, c_4, c_5 \in \mathbb{R}$ and solutions of Cauchy's functional equation $\alpha, \alpha', \beta, \beta': [0, \infty)\rightarrow \mathbb{R}$ such that
\begin{align}
\mu(P)=&c_1\,M(P)+c_2\,M([o,P])+c_3\,E([o,P])+c_4\,F([o,P])+
H_\alpha([o,P])\notag\\
+&\sum_{i=1}^{r-1} H_{\alpha'}([o,v_{i+1},v_i])+ G_\beta([o,P])+\sum_{i=1}^{r-1}G_{\beta'}([o,v_{i+1},v_i])+c_5\,\rho_{\pi/2}
\end{align}
for every $P\in \mathcal{P}^2 $ with vertices $v_1, \ldots , v_r$ visible from the origin and labeled counterclockwise.
\end{lem}

\par Theorem \ref{theo6} follows from Lemma \ref{lem3.1.1} and Lemma \ref{lem4.1.1} immediately.

\subsection{The three-dimensional case.}
\begin{lem}\label{lem4.2.1}
If $\mu: \mathcal{P}^3 \rightarrow \mathbb{M}^3$ is an $\textrm{SL}(3)$ contravariant valuation, then there exist a constant $c \in \mathbb{R}$ and a solution of Cauchy's functional equation $\xi: [0, \infty)\rightarrow \mathbb{R}$ such that
\begin{align*}
\mu (r^{\frac{1}{3}} \tilde{T}^3)=c\,I([o,r^{\frac{1}{3}} \tilde{T}^3])+2\,L_\xi([o,r^{\frac{1}{3}} \tilde{T}^3])
\end{align*}
for every $r> 0$.
\end{lem}

\emph{Proof.} \
Using  proofs  similar to Lemma \ref{lem3.2.2} and Lemma \ref{lem3.2.3}, we have $\mu(r^{\frac{1}{3}} \tilde{T}^k)=\textbf{0}$ for $k=1,2$.
Let $\phi_1, \phi_2 \in\textrm{GL}(3)$ be defined as in Definition \ref{defn2}. Combining (\ref{2.1}) and the valuation property of $\mu$ we have
\begin{align*}
\mu(r^{\frac{1}{3}} \tilde{T}^3)=\lambda^{\frac{2}{3}} {\phi_1}^{-t}\mu\left( (\lambda r)^{\frac{1}{3}} \tilde{T}^3\right) {\phi_1}^{-1}+(1-\lambda)^{\frac{2}{3}} {\phi_2}^{-t}\mu\left( ((1-\lambda)r)^{\frac{1}{3}} \tilde{T}^3\right) {\phi_2}^{-1}.
\end{align*}
Setting $\lambda=a/a+b$ and $r=a+b$ for $a, b>0$ to obtain
\begin{align*}
(a+b)^{\frac{2}{3}}\mu\left((a+b)^{\frac{1}{3}} \tilde{T}^3\right)=a^{\frac{2}{3}}{\phi_1}^{-t}\mu(a^{\frac{1}{3}} \tilde{T}^3){\phi_1}^{-1}+b^{\frac{2}{3}}{\phi_2}^{-t}\mu(b^{\frac{1}{3}} \tilde{T}^3){\phi_2}^{-1}.
\end{align*}
Next, using a similar argument as in the proof of Lemma \ref{lem3.2.4}, then there exist a constant $c \in \mathbb{R}$ and a solution of Cauchy's functional equation $\xi: [0, \infty)\rightarrow \mathbb{R}$ such that
\begin{align*}
\mu (r^{\frac{1}{3}} \tilde{T}^3)=c\,I([o,r^{\frac{1}{3}} \tilde{T}^3])+2\,L_\xi([o,r^{\frac{1}{3}} \tilde{T}^3]).
\end{align*}
\qed



\textbf{Proof of Theorem \ref{theo5}.}
It follows from Lemma \ref{lem2.2} that the expression in (\ref{1.3}) is an $\textrm{SL}(3)$ contravariant valuation. It remains to show the reverse statement.

\par First, it suffices to consider $T\in \tilde{\mathcal{T}}^{3}$ for $r > 0.$
Let $T\in \tilde{\mathcal{T}}^{3}$ with $F$ being the facet of $T$ visible from the
origin, then $[o,F]\cap T=F$ and $[o,F]\cup T=[o,T]$. Using the valuation property of $\mu$ and the inclusion-exclusion principle, we have
\begin{align}\label{4.2}
\mu([o,T])=\mu([o,F]\cup T)=\mu(T)+\mu([o,F])-\mu(F).
\end{align}
From Theorem \ref{theo2} there exist a constant $k_1 \in \mathbb{R}$ and a solution of Cauchy's functional equation $\xi_1: [0, \infty)\rightarrow \mathbb{R}$ such that
\begin{align}\label{4.3}
\mu([o,T])=k_1\, I([o,T])+2\,L_{\xi_1}([o,T])
\end{align}
and
\begin{align}\label{4.4}
\mu([o,F])=k_1\, I([o,F])+2\,L_{\xi_1}([o,F]).
\end{align}
By Lemma \ref{lem4.2.1} there exist a constant $k_2 \in \mathbb{R}$ and a solution of Cauchy's functional equation $\xi_1': [0, \infty)\rightarrow \mathbb{R}$ such that
\begin{align}\label{4.5}
\mu (F)=k_2\, I([o,F])+2\,L_{\xi_1'}([o,F]).
\end{align}
From (\ref{4.2}), (\ref{4.3}), (\ref{4.4}) and (\ref{4.5}) we obtain
\begin{align*}
\mu (T)=k_1\, I([o,T])+(k_2-k_1)\, I([o,F])+2\, L_{\xi_1}([o,T])+2\, L_{\xi_1'-\xi_1}([o,F]).
\end{align*}
Let $c_1=k_1$, $c_2=(k_2-k_1)$ and $\xi_2=\xi_1'-\xi_1$, we have
\begin{align}\label{4.6}
\mu(T)=c_1\, I([o,T])+c_2\, I([o,F])+2\, L_{\xi_1}([o,T])+2\, L_{\xi_2}([o,F]).
\end{align}

Second, let $P\in \mathcal{P}^3 \backslash \mathcal{P}_o^3,$
Triangulate $P$ into simplices $T_{1}, \cdots, T_{r}.$ Let $F_1,\ldots, F_l$ be the facets of $P$ visible from the origin.

Applying the inclusion-exclusion principle, Lemma \ref{lem 3} and (\ref{4.6}), we drive that
\begin{align*}
\mu(P)=c_1\, I([o,P])+c_2\, \sum_{i=1}^{l}I([o,F_i])+2\, L_{\xi_1}([o,P])+2\, \sum_{i=1}^{l} L_{\xi_2}([o,F_i]).
\end{align*}

Finally, for $P\in \mathcal{P}^3 $. Since it follows from Lemma \ref{lem3.2.4} that $\mu(r^{\frac{1}{3}} \tilde{T}^k)=\textbf{0}$ for $k=1,2$. Combining the triangulation, the inclusion-exclusion principle and Lemma \ref{lem4.2.1}, the assertion holds.

\qed

\subsection{The higher-dimensional case.}
In the final step, we extend Theorem \ref{theo1} to $\mathcal{P}^n$.

\begin{lem}\cite{Ma-Wang}\label{lem4.3.1}
Let $n\geq 4$. A function $\mu: \mathcal{P}^n \rightarrow \mathbb{M}_e^n$ is an $\textrm{SL}(n)$ contravariant valuation if and only if there exist solutions of Cauchy's functional equation $\zeta_1, \zeta_2:\mathbb{R} \rightarrow \mathbb{R}$ such that
\begin{equation*}
\mu(P)=L_{\zeta_1} (P)+L_{\zeta_2} ([o,P])
\end{equation*}
for every $P \in \mathcal{P}^n$, where $[o,P]$ denotes the convex hull of the origin and $P$.
\end{lem}

\par Actually, the symmetry assumption in Lemma \ref{lem4.3.1} can be removed for $n\geq 4$.

\textbf{Proof of Theorem \ref{theo4}.} \
Since $\tau_0 r T^n=r T^n$ and $\tau_0 r \tilde{T}^n=r \tilde{T}^n$, where $\tau_0 \in \textrm{SL}(n)$ comes from Lemma \ref{lem3.3.2} and $r>0$. The $\textrm{SL}(n)$ contravariance of $\mu$ implies
\begin{equation*}
\mu(r T^n)(e_i,e_j)=\mu(\tau_0 r T^n)(e_i,e_j)=\mu(r T^n)(\tau_0^{-1}e_i,\tau_0^{-1}e_j)=\mu(r T^n)(e_j,e_i)
\end{equation*}
and
\begin{equation*}
\mu(r \tilde{T}^n)(e_i,e_j)=\mu(\tau_0 r \tilde{T}^n)(e_i,e_j)=\mu(r \tilde{T}^n)(\tau_0^{-1}e_i,\tau_0^{-1}e_j)=\mu(r \tilde{T}^n)(e_j,e_i),
\end{equation*}
where $i \neq j$ and $r >0$. Hence, $\mu (r T^n)$ and $\mu (r \tilde{T}^n)$ are both symmetric.
\par Furthermore, by Lemma \ref{lem3.3.2} we obtain $\mu (r \tilde{T}^k)=\textbf{0}$ for $1\leq k \leq n-1$. Hence, the symmetry consumption of the images of $\mu$ can be removed in the proof of Ma and Wang \cite{Ma-Wang}. Apply Lemma \ref{lem4.3.1} to complete the proof.
\qed




\end{document}